\documentclass[leqno]{psp}

\usepackage{amsmath,amsfonts,amssymb,times,fancyhdr}

\newtheorem{remark}{Remark}
\newtheorem{proposition}{Proposition}
\newtheorem{lemma}{Lemma}[section]
\newtheorem{theorem}{Theorem}
\newtheorem{cor}{Corollary}

\newcommand{\be}{\begin{equation}}
\newcommand{\ee}{\end{equation}}
\newcommand{\benn}{\begin{equation*}}
\newcommand{\eenn}{\end{equation*}}

\newcommand{\li}{\text{li}}
\newcommand{\g}{\ensuremath{\gamma}}
\newcommand{\eg}{\ensuremath{e^{-\gamma}}}

\renewcommand{\b}{\ensuremath{\beta}}
\newcommand{\del}{\ensuremath{\delta}}
\newcommand{\eps}{\ensuremath{\varepsilon}}
\renewcommand{\(}{\left(}
\renewcommand{\)}{\right)}

\newcommand{\pfrac}[2]{\left(\frac{#1}{#2}\right)}

\renewcommand{\P}{\mathcal{P}}
\newcommand{\Q}{\mathcal{Q}}
\newcommand{\fl}[1]{\lfloor #1 \rfloor}
\newcommand{\hD}{\widetilde{\Delta}}

\begin{document}

\title{Generalized Euler Constants}

\author[Harold G. Diamond and Kevin Ford]{Harold G.~Diamond \nobreakand\ 
Kevin Ford\thanks{Research of the second author 
supported by National Science Foundation
grant DMS-0555367.}
 \\ University of Illinois, Urbana
}
\maketitle


\begin{abstract}
We study the distribution of a family $\{\g(\P)\}$ of generalized Euler
constants arising from integers sieved by finite sets of primes $\P$.
For $\P=\P_r$, the set of the first $r$ primes,
$\g(\P_r) \to \exp(-\g)$ as $r \to \infty$.  Calculations
suggest that $\g(\P_r)$ is monotonic in $r$, but we prove it is not.
Also, we show a connection between the distribution of $\g(\P_r) - \exp(-\g)$
and the Riemann hypothesis.
\end{abstract}


\section{Introduction} 
Euler's constant $\g =  0.5772156649 \ldots$ (also known as the
Euler-Mascheroni constant)
reflects a subtle multiplicative connection between Lebesgue measure
and the counting measure of the positive integers and appears in
many contexts in mathematics
(see e.g. the recent monograph \cite{H}).   Here we
study a class of analogues involving sieved sets of integers and
investigate some possible monotonicities.

As a first example, consider the sum of reciprocals of odd integers up
to a point $x$:  we have
\[
\sum_{\substack{n \le x \\ n \ {\rm odd}}}\,\frac 1n 
\, = \,\sum_{n\le x} \frac {1}{n} \,-\, \frac 12 
\sum_{n \le x/2}\,\frac 1n 
= \frac 12 \log x  + \frac{\g + \log 2}{2} + o(1),
\]
and we take
\[
\g_1 \!:= \lim_{x \to \infty} \Big\{\sum_{\substack{n \le x \\ 
n \ {\rm odd}}} \,\frac 1n 
-  \frac 12 \log x \Big\} = \frac{\g + \log 2}{2}\, .
\]

More generally, if ${\P}$ represents a finite set of primes, let 
\[
1_{\P}(n) \!:= \begin{cases}
1, \ \text{ if } (n, \prod_{p \in{\P}} p)=1  ,\\
0,  \ \text{ else},
\end{cases}  \ {\rm and} \quad
\delta_{\P}\!:= \lim_{x \to \infty}\frac{1}{x} \sum_{n\le x} 1_{\P}(n).
\]

A simple argument shows that $\del_{\P}= \prod_{p \in {\P}} ( 1 - 1/p )$
and that the generalized Euler constant
\[
\g(\P) \!:=  \lim_{x \to \infty} \Big\{\sum_{n \le x} 
\frac{1_{\P}(n)}{n} - \delta_{\P} \log x\Big\}
\]
exists.  We shall investigate the distribution of values of $\g(\P)$ 
for various prime sets $\P$.

\medskip

We begin by indicating two further representations of $\g(\P)$.
First, a small Abelian argument shows that it is the constant term in the
Laurent series about 1 of the Dirichlet series
\[
\sum_1^\infty 1_\P(n) n^{-s} = 
\zeta(s) \prod_{p \in \P} \big(1 - p^{-s} \big),
\]
where $\zeta$ denotes the Riemann zeta function.  That is,
\be \label{gamma-analytic}
 \g(\P) = \lim_{s \to 1} \big\{ \zeta(s) \prod_{p \in \P} 
\big(1 - p^{-s} \big) -  \frac{\delta_\P}{s-1} \big\}.
\ee

For a second representation, take $P=\prod_{p\in \P} p$.  We have
\begin{align*}
\sum_{n\le x} \frac{1_{\P}(n)}{n} &= \sum_{n\le x} \frac{1}{n} \sum_{d|(n,P)}
  \mu(d) = \sum_{d|P} \frac{\mu(d)}{d} \sum_{m\le x/d} \frac{1}{m}  \\
&= \sum_{d|P} \frac{\mu(d)}{d} \( \log (x/d) + \g + O(d/x) \) \\
&= \del_{\P} \log x - \sum_{d|P} \frac{\mu(d)\log d}{d} + \g \del_{\P} + o(1)
\end{align*}
as $x\to \infty$, where $\mu$ is the M\"obius function.  
If we apply the Dirichlet convolution identity \\ $\mu\log = 
-\Lambda * \mu$, where $\Lambda$ is the von Mangoldt function, we find that
$$
- \sum_{d|P} \frac{\mu(d)\log d}{d} = \sum_{ab|P} 
\frac{\Lambda(a) \mu(b)}{ab} 
= \sum_{p\in \P} \frac{\log p}{p} \sum_{b|P/p} \frac{\mu(b)}{b}
= \del_\P \sum_{p\in \P} \frac{\log p}{p-1}.
$$
Thus we have

\begin{proposition} \label{prop-gammaformula} 
Let $\P$ be any finite set of primes.  Then
\be\label{gammaformula}
\g(\P)  = \prod_{p \in {\P}} 
\Big( 1 - \frac {1}{p} \Big) \Big\{ \g +
\sum_{p \in {\P}} \,\frac{\log p}{p - 1} \Big\}.
\ee
\end{proposition}
We remark that this formula can also be deduced from 
\eqref{gamma-analytic} by an easy manipulation.

\medskip

It is natural to inquire about the spectrum of values 
$$
G = \{ \g(\P) \!: \P \text{ is a finite subset of primes} \}.
$$
In particular, what is
$\Gamma \!:= \inf G$ ?  
The closure of $G$ is simple to describe in terms of $\Gamma$.

\begin{proposition}  The set $G$ is dense in $[\Gamma,\infty)$.
\end{proposition}

\begin{proof}  Suppose $x>\Gamma$ and let $\P$ be a finite set of primes
with $\g(\P) < x$.  Put
$$
c = \bigl( x-\g(\P) \bigr) \prod_{p\in \P} \Big(1-\frac{1}{p}\Big)^{-1}.
$$
Let $y$ be large and let  $\P_y$ be the union of $\P$ and the primes in
$(y,e^c y]$.  By \eqref{gammaformula}, the well-known Mertens estimates,
and the prime number theorem, 
$$
\g(\P_y) = \prod_{p\in \P} \Big(1 - \frac{1}{p} \Big) \frac{\log y}{c+\log y}
\Big( 1 + O\Big(\frac{1}{\log y}\Big) \Big) \Big( \g + \sum_{p\in \P} \frac{\log p}{p-1}
+ c + O(1/\log y) \Big).
$$
Therefore, $\lim_{y\to \infty} \g(\P_y) = x$ and the proof is complete.
\end{proof}

In case ${\P}$ consists of the first $r$ primes $\{p_1,\ldots,p_r\}$, 
we replace ${\P}$ by $r$ in the preceding notation, and let $\g_r$
represent the generalized Euler constant for the integers sieved
by the first $r$ primes.  Also define $\g_0=\g(\emptyset)=\g$.
These special values play an important role in our theory of
generalized Euler constants.

\begin{table}[!h]
\caption{ \ Some Gamma Values \ (truncated)}
\begin{center}
\begin{tabular} {c c c c}
$\ \g$   = 0.57721 & & \\
$\g_1$ = 0.63518   & $\g_{11}$ = 0.56827  &  $\g_{21}$ =
0.56513  & $\g_{31}$ = 0.56385   \\  
$\g_2$ = 0.60655  &  $\g_{12}$ = 0.56783  &  $\g_{22}$ =
0.56495    & $\g_{32}$ =  0.56378  \\ 
$\g_3$ = 0.59254  &  $\g_{13}$ = 0.56745  &  $\g_{23}$ =
0.56477    & $\g_{33}$ = 0.56372   \\ 
$\g_4$ = 0.58202  &  $\g_{14}$ = 0.56694  &  $\g_{24}$ =
0.56462    & $\g_{34}$ =  0.56365  \\ 
$\g_5$ = 0.57893  &  $\g_{15}$ = 0.56649  &  $\g_{25}$ =
0.56454    & $\g_{35}$ = 0.56361   \\ 
$\g_6$ = 0.57540  &  $\g_{16}$ = 0.56619  &  $\g_{26}$ =
0.56445    & $\g_{36}$ = 0.56355   \\  
$\g_7$ = 0.57352  &  $\g_{17}$ = 0.56600  &  $\g_{27}$ =
0.56433    & $\g_{37}$ = 0.56350   \\ 
$\g_8$ = 0.57131  &  $\g_{18}$ = 0.56574  &  $\g_{28}$ =
0.56420    & $\g_{38}$ = 0.56345    \\  
$\g_9$ = 0.56978  &  $\g_{19}$ = 0.56555  &  $\g_{29}$ =
0.56406    & $\g_{39}$ = 0.56341   \\ 
$\!\!\g_{10}$ = 0.56913 & $\g_{20}$ = 0.56537  & $\g_{30}$ =
0.56391    & $\g_{40}$ = 0.56336   \\ 
& & &  \!\!$e^{-\g}$ = 0.5614594835\ldots \\ 
\end{tabular}
\end{center}
\end{table}

The next result will be proved in Section \ref{sec:comparison}.

\begin{theorem}\label{comparison}
Let $\P$ be a finite set of primes.  For some $r$, $0\le r
\le \# \P$, we have $\g(\P) \ge \g_r$.  Consequently, 
$\Gamma = \inf_{r\ge 0} \gamma_r$.
\end{theorem}

Applying Mertens' well-known formulas for sums and products of primes
to \eqref{gammaformula}, we find that
\be \label{gamma-limit}
\g_r \sim \frac{\eg}{\log p_r}\, \{ \log p_r + O(1) \}\,
\sim \, \eg  \quad (r \to \infty) \, .
\ee
In particular, $\Gamma \le \eg$ and $G$ is dense in $[\eg,\infty)$.

Values of $\g_r$ for all $r$ with $p_r \le 10^9$ were computed to high
precision using PARI/GP.  For all such $r$, $\g_r > \eg$ and
$\g_{r+1} < \g_r$.  It is natural to ask if these trends persist.
That is, (1) Is the sequence of $\g_r$'s indeed decreasing for all $r \ge 1$? 
(2) If the $\g_r$'s
oscillate, are any of them smaller than $e^{-\g}$, i.e., is
$\Gamma < \eg$ ? \emph{
We shall show that the answer to (1) is No and the answer to (2) is No
or Yes depending on whether the Riemann Hypothesis (RH) is
true or false.}

\begin{theorem}\label{oscillation}
There are infinitely many integers $r$ with $\g_{r+1} > \g_r$
and infinitely many integers $r$ with $\g_{r+1} < \g_r$.
\end{theorem}

Theorem \ref{oscillation} will be proved in Section \ref{osc}.
There we also argue that the smallest $r$ satisfying $\g_{r+1}>\g_r$
is probably larger than $10^{215}$, and hence no amount of
computer calculation (today) would detect this phenomenon.  This behavior is
closely linked to the classical problem of locating 
sign changes of $\pi(x)-\li(x)$,
where $\pi(x)$ is the number of primes $\le x$ and 
$$
\li(x) = \int_0^x \frac{dt}{\log t} = \lim_{\eps\to 0^+}
\( \int_0^{1-\eps} + \int_{1+\eps}^x \) \frac{dt}{\log t}
$$
is Gauss' approximation to $\pi(x)$.

Despite the oscillations, $\{\g_r\}$ can be shown (on RH) to
approach $\eg$ from above.  If RH is false, $\{\g_r\}$ assumes values
above and below $\eg$ (while converging to this value).

\begin{theorem}\label{onRH}
Assume RH.  Then $\g_r > \eg$ for all $r\ge 0$.  Moreover, we have
\be\label{grasym}
\g_r = \eg \( 1 + \frac{g(p_r)}{\sqrt{p_r} (\log p_r)^2} \),
\ee
where $1.95 \le g(x) \le 2.05$ for large $x$.
\end{theorem}

As we shall see later,
$\limsup_{x\to\infty} g(x) > 2$ and $\liminf_{x\to\infty} g(x) < 2$.

\begin{theorem}\label{RHfalse}
Assume RH is false.  Then $\g_r < \eg$ for infinitely many $r$.
In particular, $\Gamma < \eg$.
\end{theorem}

\begin{cor}
The Riemann Hypothesis is equivalent to the statement ``$\g_{r} > \eg$ for all
$r\ge 0$.''
\end{cor}

It is relatively easy to find reasonable, unconditional 
lower bounds on $\Gamma$ by making use of Theorem \ref{comparison},
Proposition \ref{prop-gammaformula}
and explicit bounds for counting functions of primes.
By Theorem 7 of \cite{RS}, we have
$$
\prod_{p\le x} \( 1 - \frac{1}{p} \) > \frac{e^{-\g}}{\log x} \( 1 - 
\frac{1}{2\log^2 x} \) \qquad (x\ge 285).
$$
Theorem 6 of \cite{RS} states that
$$
\sum_{p\le x} \frac{\log p}{p} > \log x - \g - \sum_p \frac{\log p}{p(p-1)}
- \frac{1}{2\log x} \qquad (x>1).
$$
Using Proposition \ref{prop-gammaformula} and writing $\frac{1}{p-1} =
\frac{1}{p} + \frac{1}{p(p-1)}$, we obtain for $x=p_r\ge 285$ the bound
\begin{align*}
\g_r &\ge \frac{e^{-\g}}{\log x} \( 1 - \frac{1}{2\log^2 x} \) \(
\g + \sum_{p\le x} \frac{\log p}{p} + \sum_p \frac{\log p}{p(p-1)} - 
\sum_{p \ge x+1} \frac{\log p}{p(p-1)} \) \\
&\ge \frac{e^{-\g}}{\log x} \( 1 - \frac{1}{2\log^2 x} \) \( \log x - \frac{1}
{2\log x} - \frac{(x+1)(1+\log x)}{x^2} \).
\end{align*}
In the last step we used
$$
\sum_{p\ge x+1} \frac{\log p}{p(p-1)} < \frac{x+1}{x} \int_x^\infty 
\frac{\log t}{t^2}\, dt = \frac{(x+1)(1+\log x)}{x^2}.
$$

By the aforementioned computer calculations, $\g_r > e^{-\g}$ when
$p_r < 10^9$, and for $p_r > 10^9$ the bound given above implies that
$\g_r \ge 0.56$.  Therefore, we have unconditionally
$$
\Gamma \ge 0.56.
$$

Better lower bounds can be achieved by utilizing longer computer
calculations, better bounds for prime counts \cite{RR}, and some of
the results from \S \ref{infRH} below, especially \eqref{gammaRHD2}.

%
\section{An extremal property of $\boldsymbol{\{\g_r\}}$}
\label{sec:comparison}
%

In this section we prove Theorem \ref{comparison}.  Starting with an
arbitrary finite set $\P$ of primes, we perform a sequence of
operations on $\P$, at each step either removing the largest prime
from our set or replacing the largest prime with a smaller one.  We
stop when the resulting set is the first $r$ primes, with $0\le r\le
\# \P$.  We make strategic choices of the operations to create a
sequence of sets of primes $\P_0=\P, \P_1, \ldots, \P_k$, where
$$
\g(\P_0) > \g(\P_1) > \cdots > \g(\P_k)
$$
with $\P_k = \{p_1, p_2,\ldots, p_r\}$, the first $r$ primes.

The method is simple to describe.  Let $\Q=\P_j$, which is not equal
to any set $\{p_1, p_2, \ldots,p_s\}$, and with largest element $t$.
Let $\Q' = \Q \backslash \{t\}$.  If $\g(\Q') < \g(\Q)$, we set
$\P_{j+1}=\Q'$.  Otherwise, we set $\P_{j+1} = \Q' \cup \{u\}$, where
$u$ is the smallest prime not in $\Q$.  We have $u<t$ by assumption.  
It remains to show in the latter case that 
\be\label{QQ''}
\g(\P_{j+1}) < \g(\Q).
\ee
By \eqref{gammaformula}, for any prime $v \not\in \Q'$,
\be\label{addprime}
\g(\Q' \cup \{v\}) = \g(\Q') \( 1 - \frac{1}{v} + \frac{\log v}{vA} \) =:
\!\g(\Q') f(v), \quad
A\!: = \g + \sum_{p\in\Q'} \frac{\log p}{p-1} \,.
\ee
Observe that $f(v)$ is strictly increasing for $v<e^{A+1}$ and strictly
decreasing for $v>e^{A+1}$, and $\lim_{v\to\infty} f(v)=1$.  Thus
$f(v) > 1$ for $e^{A+1} \le v < \infty$.  Since
$\g(\Q)=\g(\Q') f(t) \le \g(\Q')$, we have $f(t) \le 1$.
It follows that $u<t \le e^{A+1}$ and hence $f(u) < f(t) \le 1$.
Another application of \eqref{addprime}, this time with $v=u$,
proves \eqref{QQ''} and the theorem follows.

%
\section{The $\boldsymbol{\g_r}$'s \ are not monotone}\label{osc}
%

Define
\be\label{Ax}
A(x) \!:= \g + \sum_{p\le x} \frac{\log p}{p-1}.
\ee
By \eqref{addprime}, we have
$$
\g_{r+1} = \g_r \( 1 - \frac{1}{p_{r+1}} + \frac{\log p_{r+1}}{p_{r+1}
 A(p_{r})} \),
$$
thus
\be\label{monotcond}
\g_{r+1} \le \g_r \quad \iff \quad A(p_{r}) \ge \log p_{r+1}.
\ee

\begin{theorem}\label{Axosc}
We have 
$A(x)-\log x = \Omega_\pm (x^{-1/2} \log\log\log x)$.
\end{theorem}

\begin{proof}
First introduce
\be\label{Deltax}
\Delta(x) \!:= \sum_{p \le x} \,\frac{\log p}{p - 1} 
- \sum_{n \le x}  \,\frac{\Lambda(n)}{n} \quad  {\rm and} \quad
\theta(x) \!:= \sum_{p \le x} \log p \, .
\ee
Then
\be\label{Delta1}
\begin{split}
\Delta(x)&= \sum_{p \le x} \, (\log p) \sum_{\alpha \ge 1} \frac{1}{p^\alpha}
-  \sum_{p^\alpha \le x} \, (\log p) \frac{1}{p^\alpha}  \\
&= \sum_{p \le x} \,  (\log p)  \sum_{\alpha \ge \fl{\log x /\log p} +1} 
\frac{1}{p^\alpha} = 
\sum_{p \le x} \,  \frac{\log p}{p-1} \, p^{-\fl{\log x / \log p}} \ge 0.
\end{split}
\ee
Since $p^{\fl{\log x / \log p}} \ge x/p$, we have
\be\label{Delta2}
\Delta(x) = \sum_{\sqrt{x} < p \le x} \frac{\log p}{p(p-1)} +
\sum_{x^{1/3} < p \le x^{1/2}} \frac{\log p}{p^2(p-1)} +
O \bigg( \sum_{p\le x^{1/3}} \frac{\log p}{x} \bigg).
\ee
Aside from an error of $O(x^{-1})$, the first sum is
$$
\sum_{p>\sqrt{x}} \frac{\log p}{p^2} = - \frac{\theta(\sqrt{x})}{x} +
\int_{\sqrt{x}}^\infty \frac{2\theta(t)}{t^3}\, dt 
= x^{-1/2} + O(x^{-1/2}\log^{-3} x),
$$
using the bound $|\theta(x)-x| \ll x\log^{-3} x$ which follows from
the prime number theorem with a suitable error term.
The second sum and error term in \eqref{Delta2} are each $O(x^{-2/3})$,
and we deduce that
\be\label{Deltaasym}
\Delta(x) = \frac{1}{\sqrt{x}} + O\pfrac{1}{\sqrt{x} (\log x)^3}.
\ee

\begin{remark}
Assuming RH and using the von Koch bound $|\theta(x)-x| \ll
\sqrt{x} \log^2 x$, we obtain the sharper estimate $\Delta(x)=x^{-1/2}
+ O(x^{-2/3})$.
\end{remark}

By \eqref{Deltaasym},
\be\label{Ax1}
A(x) = \g  + \sum_{n\le x} \frac{\Lambda(n)}{n} + \frac{1}{\sqrt{x}} 
+ O\pfrac{1}{\sqrt{x} (\log x)^3}.
\ee

To analyze the above sum, introduce
\be \label{R-def}
R(x) \!:= \sum_{n\le x} \frac{\Lambda(n)}{n} - \log x + \g.
\ee
For $\Re s > 0$, we compute the Mellin transform
\be\label{MTR}
\int_{1}^\infty x^{-s-1} R(x)\, dx = -\frac{1}{s} \, \frac{\zeta'}{\zeta}(s+1)
- \frac{1}{s^2} + \frac{\g}{s}.
\ee
The largest real singularity of the function on the right comes from the
trivial zero of $\zeta(s+1)$ at $s=-3.$

Let $\rho = \b+i\tau$ represent a generic nontrivial zero of $\zeta(s)$ 
-- we avoid use of $\g$ for $\Im \rho$ for obvious reasons.  
If RH is false, there is a zero $\b+i\tau$ of 
$\zeta$ with $\b>1/2$, and a straightforward application of
Landau's Oscillation Theorem (\cite{BD}, Theorem 6.31) gives
$R(x) = \Omega_{\pm} (x^{\b-1-\eps})$
for every $\eps>0$.   In this case,  $A(x) - \log x =
\Omega_{\pm} (x^{\b-1-\eps})$, which is stronger than the assertion
of the theorem.

If RH is true, we may analyze $R(x)$ via the explicit formula
\be
\label{Rexplicit} 
R_0(x) \!:= \frac 12 \{R(x^-) + R(x^+)\} =-\sum_\rho 
\frac{x^{\rho-1}}{\rho-1} + \sum_{n=1}^\infty\,\frac{1}{2n+1} \, x^{-2n-1} \,, 
\ee 
where $\sum_\rho$ means $\lim_{T\to \infty} \sum_{|\rho| \le T}$.  Equation
\eqref{Rexplicit} is deduced in a standard way from \eqref{MTR} by
contour integration, and $\lim_{T\to\infty} \sum_{|\rho| \le T}$ converges
boundedly for $x$ in any (fixed) compact set contained in $(1, \,
\infty)$. (cf.~\cite{Dav}, Ch.~17, where a similar formula is given
for $\psi_0(x)$, as we now describe.)

In showing that
\[
\psi(x) \!:= \sum_{n \le x} \Lambda(n) = x +
\Omega_\pm(x^{1/2}\log \log \log x),
\]
Littlewood \cite{Li} (also cf.~\cite{Dav}, Ch.~17) used the
analogous explicit formula 
\[
\psi_0(x) \!:=
\frac 12 \{\psi(x^+) + \psi(x^-)\} = x  -  \frac{\zeta'}{\zeta}(0)
 - \sum_\rho \frac{x^\rho}{\rho} + \sum_{n=1}^\infty\, 
\frac{1}{2n} \, x^{-2n} 
\]
and proved that
\[
\sum_\rho \frac{x^{\rho}}{\rho} 
=\Omega_\pm(\sqrt{x} \log \log \log x).
\]
Forming a difference of normalized sums over the non-trivial zeros $\rho$,
we obtain
\[
\Big| \sum_\rho \frac{x^{\rho-1/2}}{\rho} - 
 \sum_\rho \frac{x^{\rho-1/2}}{\rho-1} \Big| 
\ll \sum_\rho\Big| \frac{1}{\rho (\rho-1)} \Big| \ll 1.
\]
Thus
\[
\sum_\rho \frac{x^{\rho-1}}{\rho-1} = x^{-1/2} \,
\Omega_\pm(\log \log \log x),
\]
and hence $R(x) = \Omega_\pm(x^{-1/2} \log \log \log x)$.

Therefore, in both cases (RH false, RH true), we have
$$
R(x) = \Omega_\pm(x^{-1/2} \log \log \log x),
$$
and the theorem follows from \eqref{Ax1} and \eqref{R-def}.
\end{proof}

By Theorem \ref{Axosc}, there are arbitrarily large values of $x$
for which $A(x) < \log x$.  
If $p_{r}$ is the largest prime $\le x$, then 
$$
A(p_r) = A(x) < \log x < \log p_{r+1}
$$
for such $x$.  This
implies by \eqref{monotcond} that $\g_{r+1} > \g_r$.
For the second part of Theorem \ref{oscillation}, take $x$ large and satisfying
$A(x) > \log x + x^{-1/2}$ and let $p_{r+1}$ be the largest prime 
$\le x$.  By Bertrand's postulate, $p_{r+1} \ge x/2$.  Hence
\[
A(p_r) = A(p_{r+1})  -\frac{\log p_{r+1}}{p_{r+1}-1}  
\ge \log x + x^{-1/2}  -\frac{\log x}{x/2-1}
> \log x \ge \log p_{r+1},
\]
which implies $\g_{r+1} < \g_r$.
\bigskip

Computations with PARI/GP reveal that $\g_{r+1} < \g_r$ 
for all $r$ with $p_r < 10^9$.  
By \eqref{Ax1}, to find $r$ such that $\g_{r+1}>\g_r$, we need 
to search for values of $x$ essentially satisfying $R(x) < -x^{-1/2}$.
By \eqref{Rexplicit}, this boils down to finding values of $u=\log x$
such that
$$
\sum_{\rho=\b+i\tau} \frac{e^{iu\tau}}{i\tau-1/2} > 1.
$$ 
Of course, the smallest zeros of $\zeta(s)$ make the greatest
contributions to this sum.

Let $\ell(u)$ be the truncated version of the preceding sum taken over
the zeros $\rho$ with $|\Im \rho| \le T_0 := 1132490.66$ (approximately 2
million zeros with positive imaginary part, together with their
conjugates).  A table of these zeros, accurate to within $3\cdot 10^{-9}$,
is provided on Andrew Odlyzko's web page

{\tt http://www.dtc.umn.edu/}$\sim${\tt odlyzko/zeta\_tables/index.html}.

In computations of $\ell(u)$, the errors in the values of the zeros
contribute a total error of at most
$$
(3\cdot 10^{-9}) u \sum_{|\Im \rho| \le T_0} \frac{1}{|\rho|} \le
(4.5\cdot 10^{-7})u.
$$
Computation using $u$-values at increments of $10^{-5}$ and an
early abort strategy for $u$'s having too small a sum over the first
1000 zeros, indicates that $\ell(u) \le 0.92$ for $10 \le u\le 495.7$.
Thus, it seems likely that the first $r$ with
$\g_{r+1}>\g_r$ occurs when $p_r$ is of size at least $e^{495.7}
\approx 1.9 \times 10^{215}$.  There is a possibility that the first
occurence of $\g_{r+1}>\g_r$ happens nearby, as
$\ell(495.702808) > 0.996$.
Going out a bit further, we find that $\ell(1859.129184) > 1.05$, and
an averaging method of R.~S.~Lehman \cite{Leh} can be used to prove
that $\g_{r+1}>\g_r$ for many values of $r$ in the vicinity of
$e^{1859.129184} \approx 2.567 \times 10^{807}$.  Incidentally, for
the problem of locating sign changes of $\pi(x)-\li(x)$, one must find
values of $u$ for which (essentially)
$$
\sum_{\rho=\b+i\tau} \frac{e^{iu\tau}}{i\tau+1/2} < -1.
$$
A similarly truncated sum over zeros with $|\Im \rho| \le
600,000$ first attains values less than $-1$ for positive $u$
values when $u \approx 1.398 \times 10^{316}$ \cite{BH}.

%
\section{Proof of Theorem \ref{onRH}}\label{infRH}
%

Showing that $\g_r > \eg$ 
for \emph{all} $r\ge 0$ under RH requires explicit estimates for prime
numbers.  Although sharper estimates are known (cf. \cite{RR}),
older results of Rosser and Schoenfeld suffice for our purposes.  
The next lemma follows from Theorems 9 and 10 of \cite{RS}. 

\begin{lemma}\label{RStheta}
We have 
$\,\theta(x) \le 1.017 x\,$ for $\,x>0\,$ and $\, \theta(x) \ge 0.945 \,x\,$
for $\,x\ge 1000$.
\end{lemma}

The preceeding lemma is unconditional.  On RH, we can do better for large $x$,
such as the following results of Schoenfeld (\cite{S76}, 
Theorem 10 and Corollary 2).

\begin{lemma}\label{Sch}
Assume RH.  Then
$$
|\theta(x)-x| < \frac{\sqrt{x} \log^2 x}{8\pi} \qquad (x\ge 599)
$$
and
$$
|R(x)| \le \frac{3\log^2 x + 6\log x + 12}{8\pi \sqrt{x}} \qquad (x\ge 8.4).
$$
\end{lemma}

Mertens' formula in the form
\[
- \sum_{p\le x} \log(1-1/p) = \log \log x + \g + o(1)
\]
and a familiar small calculation give
\[
- \sum_{p\le x} \log(1-1/p) - \sum_{n\le x} \frac{\Lambda(n)}{n\log n}
= \sum_{\substack{p\le x \\ p^a > x}} \frac{1}{ap^a}
= O\Big(\frac{1}{\log x}\Big) = o(1).
\]
It follows that
\be \label{lam-over-log}
\sum_{n\le x} \frac{\Lambda(n)}{n\log n} = \log \log x + \g + o(1).
\ee

We can obtain an exact expression for the last sum in terms of $R$
(defined in \eqref{R-def}) by integrating by parts:
\begin{align*}
\sum_{n\le x} \frac{\Lambda(n)}{n\log n} 
&= \int_{2^-}^x \frac{dt/t + dR(t)}{\log t}  \\
&= \log \log x - \log \log 2 + \frac{R(x)}{\log x} - \frac{R(2)}{\log 2}
+ \int_2^x \frac{R(t) \, dt}{t \log^2 t} \\
&=\log \log x + c
+ \frac{R(x)}{\log x} - \int_x^\infty \frac{R(t) \, dt}{t \log^2 t} \,, 
\end{align*}
where
\[
c \!:= \int_2^\infty \frac{R(t) \, dt}{t \log^2 t}
- \frac{R(2)}{\log 2} - \log \log 2 = \g,
\] 
by reference to \eqref{lam-over-log} and the relation $R(x)=o(1)$.  Thus
\be\label{logzeta}
\sum_{n\le x} \frac{\Lambda(n)}{n\log n} = \log\log x + \g + \frac{R(x)}{\log
  x} - \int_{x}^\infty \frac{R(t)}{t\log^2 t}\, dt.
\ee

Let $H(x)$ denote the integral in \eqref{logzeta} and define
\be\label{hDx}
\hD(x) \!:= \sum_{\substack{p\le x \\ p^a > x}} \frac{1}{ap^a}.
\ee

Using \eqref{gammaformula}, engaging nearly all the preceding 
notation and writing $p_r = x$, we have
\be\label{gammaRHD}
\g_r = \frac{\eg}{\log x} \exp \Big\{ - \frac{R(x)}{\log x} + H(x) - \hD(x)
\Big\} \bigl( \log x + R(x) + \Delta(x) \bigr).
\ee
We use Lemmas \ref{RStheta} and \ref{Sch} to obtain explicit estimates
for $H(x)$, $\,\Delta(x)$, and $\,\hD(x)$.  

We shall show below that  $\Delta(x) - \hD(x)\log x \ge 0$. 
It is crucial for our arguments that this difference be small. Also,
although one may use Lemma \ref{Sch} to bound $H(x)$, we shall obtain
a much better inequality by using the explicit formula \eqref{Rexplicit}
for $R_0$ (which agrees with $R$ a.e.).

\begin{lemma}\label{Hlem}
Assume RH.  Then
$$
|H(x)| \le \frac{0.0462}{\sqrt{x} \log^2 x} \Big(1 + \frac{4}{\log x}
\Big) \qquad (x\ge 100).
$$
\end{lemma}

\begin{proof}
Since $R(x)=o(1)$ by the Prime Number Theorem, we see
that the integral defining $H$ converges absolutely.  We write
\[
H(x) = \lim_{X\to \infty} \int_x^X \frac{R(t)}{t \log^2 t} \, dt,
\]
and treat the integral for $H$ as a finite integral in justifying
term-wise operations.

We now apply the explicit formula \eqref{Rexplicit} for $R$.  For $t\ge 100$, 
$$
\sum_{n=1}^\infty \frac{t^{-2n-1}}{2n+1} \le \frac{0.34}{t^3}
$$
and thus
\be\label{Htrvialzeros}
\int_x^\infty \frac{1}{t\log^2 t} \sum_{n=1}^\infty \frac{t^{-2n-1}}{2n+1}
\, dt 
< \frac{0.12}{x^3\log^2 x} \le \frac{0.0000012}{x^{1/2}\log^2 x}.
\ee

The series over zeta zeros in \eqref{Rexplicit} converges boundedly to
$R_0(x)$ as $T \to \infty $ for $x$ in a compact region; by the
preceding remark on the integral defining $H$, we can integrate the
series term-wise.  For each nontrivial zero $\rho$, integration by parts gives
\[
\int_x^\infty \frac{t^{\rho-2}}{\log^2 t}\, dt 
=\frac{-x^{\rho-1}}{(\rho-1)\log^2 x} + \frac{2}{\rho-1} \int_x^\infty
\frac{t^{\rho-2}}{\log^3 t}\, dt
\]
and thus

\begin{align*}
\left| \int_x^\infty \frac{t^{\rho-2}}{\log^2 t}\, dt \right|
&\le \frac{x^{-1/2}}{|\rho-1| \log^2 x} + \frac{2}{|\rho-1|\log^3 x}
\int_x^\infty t^{-3/2}\, dt \\
&\le \frac{x^{-1/2}}{|\rho-1|\log^2 x} \(1 + \frac{4}{\log x} \).
\end{align*}
Since RH is assumed true, we have by \cite{Dav}, Ch. 12, (10) and (11),
$$
\sum_{\rho} \frac{1}{|\rho-1|^2} = \sum_{\rho} \frac{1}{|\rho|^2}
= 2 \sum_{\rho} \frac{\Re \rho}{|\rho|^2} = 
2 + \g - \log 4\pi  =  0.0461914 \ldots.
$$
Putting these pieces together, we conclude that
\begin{align*}  
|H(x)| &\le \sum_{\rho} \frac{1}{|\rho-1|} \left| \int_x^\infty 
\frac{t^{\rho-2}}{\log^2 t}\, dt \right| +
\frac{0.0000012}{\sqrt x \, \log^2 x}  \\
&\le \frac{0.0461915+0.0000012}{\sqrt{x}\log^2 x} \( 1 + \frac{4}{\log x}
\). 
\end{align*}
\end{proof}

Under assumption of RH, we have
\begin{equation}  \label{H-formula}
H(x) = \frac{x^{-1/2}}{\log^2 x} \Big\{\sum_\rho
\frac{x^{i\tau}}{(\rho-1)^2}+\frac{4 \vartheta}{\log x}\sum_\rho
\frac{1}{|\rho-1|^2} + \frac{0.12 \vartheta'}{x^{5/2}} \Big\}\,,
\end{equation}
where $|\vartheta|\le 1$ and $|\vartheta'|\le 1$.
The series are each absolutely summable, and so the first series is an
almost periodic function of $\log x$. Thus the values this series
assumes are (nearly) repeated infinitely often.  The other two terms
in \eqref{H-formula} converge to 0 as $x \to \infty$.  Also, the mean
value of $H(x) x^{1/2} \log^2 x$ is 0 (integrate the first series); thus
the first series in \eqref{H-formula} assumes both positive and
negative values.  The $\limsup$ and $\liminf$ of $H(x) x^{1/2} \log x$ are
equal to the $\limsup$ and $\liminf$ of the first series in
\eqref{H-formula}, and we have
\[
H(x)=\Omega_\pm(x^{-1/2} (\log x)^{-2}).
\]

If one assumes that the zeros $\rho$ in the
upper half-plane have imaginary parts which are linearly independent
over the rationals (unproved even under RH, but widely believed), 
then Kronecker's theorem implies that
$$
\limsup_{x\to \infty} H(x) \sqrt{x} (\log x)^2 = 2+\g-\log 4\pi, \quad
\liminf_{x\to \infty} H(x) \sqrt{x} (\log x)^2 = -(2+\g-\log 4\pi).
$$

Continuing to assume RH but making no linear independence assumption
on the $\tau$'s, we can show that 
\[
\liminf_{x\to \infty} H(x) \sqrt{x} (\log x)^2 \le
- \sum_\rho |\rho-1|^{-2} + \frac 12 \sum_\rho |\rho-1|^{-4} < -0.04615,
\]
which is close to $-(2+\g-\log 4\pi)$.  Indeed, for $x=1$, the first
series in \eqref{H-formula} equals $\sum_{\rho} (\rho-1)^{-2}$, and
by almost periodicity this value is nearly repeated infinitely often.
Also,
\[
\frac{1}{(\rho-1)^2} + \frac{1}{(\overline \rho-1)^2} 
+ \frac{2}{|\rho - 1|^2} = \frac{1}{|\rho - 1|^4},
\]
so that
\[
\sum_{\rho} \frac{1}{(\rho-1)^2} = - \sum_{\rho} \frac{1}{|\rho-1|^2} +
 \sum_{\rho} \frac{1/2}{|\rho-1|^4}.
\]

The next two lemmas are unconditional; i.e. they do not depend on RH.
We do not try to obtain the sharpest estimates here.

\begin{lemma}\label{Deltalem}
We have 
$$
\Delta(x) \le \frac{3.05}{\sqrt{x}} \qquad (x\ge 10^6).
$$
\end{lemma}

\begin{proof}
Using \eqref{Delta1} and the upper bound for $\theta(x)$ given in 
Lemma \ref{RStheta},
\begin{align*}
\Delta(x) &\le \frac{\sqrt{x}}{\sqrt{x}-1} \sum_{p>\sqrt{x}} \frac{\log
  p}{p^2} + 2 \sum_{p\le \sqrt{x}} \frac{\log p}{x} \\
&\le \frac{1000}{999} \( - \frac{\theta(\sqrt{x})}{x} + \int_{\sqrt{x}}^\infty
  \frac{2\theta(t)}{t^3}\, dt \) + \frac{2\theta(\sqrt{x})}{x} \\
&\le 1.017 \( 4 -\frac{1000}{999}\) x^{-1/2} < 3.05 x^{-1/2}.
\end{align*}
\end{proof}

\begin{lemma}\label{DeltahD}
We have 
\begin{align}
& \Delta(x)/\log x \ge \hD(x) \quad \ (x>1)  \label{DeltahD1}  \\
&\frac{\Delta(x)}{\log x} - \hD(x) = \frac{2}{\sqrt{x}\log^2 x} +
O\pfrac{1}{\sqrt{x}\log^3 x} \quad \ (x\ge 2)  \label{DeltahD2}  \\
&\frac{\Delta(x)}{\log x} - \hD(x) \ge \frac{1.23}{\sqrt{x}\log^2 x} 
\quad \ (x\ge 10^6). \label{DeltahD3}
\end{align}
\end{lemma}

\begin{proof}
By \eqref{Delta1}, we have
\begin{equation}  \label{Del-diff}
\frac{\Delta(x)}{\log x} - \hD(x) = \sum_{p\le x} \sum_{a >
\frac{\log x}{\log  p}} \frac{1}{p^a} \( \frac{\log p}{\log x} -
\frac{1}{a} \). 
\end{equation}
Each summand on the right side is clearly positive, proving the first 
part of the lemma.

As shown in the proof of \eqref{Delta2}, the summands of $\Delta(x)$
associated with exponents $a\ge 3$ make a total contribution of
$O(x^{-2/3})$. Thus the corresponding summands in \eqref{Del-diff}
contribute $O(x^{-2/3}/\log x)$.  We handle the remaining term by
partial summation, writing
\begin{align} \label{DhD1}
\sum_{\sqrt{x}<p\le x} \!p^{-2} \( \frac{\log p}{\log x} - \frac12\) 
&= \int_{\sqrt x}^x \frac{1}{t^2}\Big\{ \frac{1}{\log x} - \frac{1}{2
  \log t} \Big\} d\theta(t) \\ 
&=\frac{\theta(x)}{2x^2\log x} + \int_{\sqrt{x}}^x \theta(t)t^{-3} 
\Big\{ \frac{2}{\log x} - \frac{1}{\log t}  - \frac{1}{2\log^2 t}
\Big\} dt.  \notag
\end{align}

Using the prime number theorem with an error term $\theta(t) - t \ll
t \log^{-2} t$, the left side of \eqref{DhD1} is seen to be 
\begin{align*}
&=  O\pfrac{1}{\sqrt{x}\log^3 x} + \int_{\sqrt{x}}^\infty \Big\{ \Big(
  \frac{2}{t^2\log x} - \frac{1}{t^2\log t} - \frac{1}{t^2\log^2 t} \Big)
  + \frac{1}{2t^2\log^2 t} \Big\}\, dt \\
&=  O\pfrac{1}{\sqrt{x}\log^3 x} +  \int_{\sqrt{x}}^\infty 
  \frac{1}{2t^2\log^2 t}\, dt \\
&= \frac{2}{\sqrt{x} \log^2 x} + O\pfrac{1}{\sqrt{x}\log^3 x}.
\end{align*}
proving the second part of the lemma.

The proof of \eqref{DeltahD1} shows the expression in \eqref{DhD1} is
a valid lower bound for $\Delta(x)/\log x - \hD(x)$.  Inserting the
estimates from Lemma \ref{RStheta} and applying integration by parts
gives, for $x\ge 10^6$,
\begin{align*}
\frac{\Delta(x)}{\log x} - \hD(x) &\ge \frac{0.4725}{x\log x}
  + 0.945 \int_{\sqrt{x}}^x \frac{2/\log x - 1/\log t}{t^2}\, dt - 
  \frac{1.017}{2} \int_{\sqrt{x}}^x \frac{dt}{t^2\log^2 t} \\
&= - \frac{0.4725}{x\log x} + 0.4365 \int_{\sqrt{x}}^x \frac{dt}{t^2\log^2 t}.
\end{align*}
Another application of integration by parts yields
\begin{align*}
\int_{\sqrt{x}}^x \frac{dt}{t^2\log^2 t} &= \frac{4}{\sqrt{x}\log^2 x}
- \frac{1}{x\log^2 x} -\int_{\sqrt{x}}^x \frac{2 dt}{t^2 \log^3 t} \\
&\ge \frac{4}{\sqrt{x}\log^2 x} - \frac{1}{x\log^2 x} 
- \frac{16}{\sqrt{x} \log^3 x} \\
&\ge \frac{2.84}{\sqrt{x}\log^2 x}.
\end{align*}
Finally,
$$
\frac{1}{x\log x} \le \frac{\log 10^6}{1000} \frac{1}{\sqrt{x}\log^2 x}.
$$
Combining the estimates, we obtain the third part of the lemma. 
\end{proof}

We are now set to complete the proof of Theorem \ref{onRH}.
A short calculation using PARI/GP verifies that $\g_r > \eg$ for
$p_r < 10^6$.  Assume now that $x=p_r \ge 10^6$.  By \eqref{gammaRHD},
\be\label{gammaRHD2}
\g_r = \eg \(1 + \frac{R(x)+\Delta(x)}{\log x} \) \exp \left\{
- \frac{R(x)+\Delta(x)}{\log x} \right\} \exp \left\{ 
\frac{\Delta(x)}{\log x} - \hD(x) + H(x) \right\}.
\ee
By Lemmas \ref{Sch} and \ref{Deltalem},
$$
\frac{|R(x)|+\Delta(x)}{\log x} \le \frac{3\log^2 x+6\log x+12+24.4\pi}
{8\pi \sqrt{x} \log x} \le \frac{0.1556\log x}{\sqrt{x}} \le 0.00215.
$$
By Taylor's theorem applied to $-y+\log(1+y)$, if $|y| \le 0.00215$ then
$e^{-y}(1+y) \ge e^{-0.501y^2}$.  This, together with Lemmas \ref{Hlem}
and \ref{DeltahD}, yields 
$$
\g_r \ge \eg \exp \left\{ -0.01213 \frac{\log^2 x}{x} + \frac{1.17}
{\sqrt{x}\log^2 x} \right\}.
$$
Since $x^{-1/2} \log^4 x$ is decreasing for $x\ge e^8$, 
$$
\frac{\log^2 x}{x} \le \frac{\log^4 10^6}{1000} \frac{1}{\sqrt{x}\log^2 x}
\le \frac{36.431}{\sqrt{x}\log^2 x}.
$$
We conclude that
$$
\g_r \ge \eg \exp \left\{ \frac{0.728}{\sqrt{x}\log^2 x} \right\} \qquad
(x\ge 10^6),
$$
which completes the proof of the first assertion.

By combining \eqref{gammaRHD2} with Lemma \ref{Sch}, Lemma \ref{Deltalem} and
\eqref{DeltahD2}, we have
$$
\g_r = e^{-\g} \exp \left\{ H(x) + \frac{2}{\sqrt{x} \log^2 x} + O\pfrac{1}
{\sqrt{x} \log^3 x} \right\}.
$$
Lemma \ref{Hlem} implies that 
$$
|H(x)| \le \frac{0.047}{\sqrt{x} \log^2 x}
$$
for large $x$, 
and this proves \eqref{grasym}.  By the commentary following the proof of Lemma
\ref{Hlem}, we see that $\liminf g(x) < 2$ and $\limsup g(x) > 2$.

%
\section{Analysis of $\boldsymbol{\g_r}$ if RH is false}\label{MTH}
%

Start with \eqref{gammaRHD2} and note that $e^{-y} (1+y) \le 1$.
Inserting the estimates from Lemma \ref{DeltahD}
gives
\be\label{grbigosc}
\g_r \le \eg \exp \Big\{H(x) + O\Big(\frac{1}{\sqrt{x}\log^2 x}\Big) \Big\}.
\ee
Our goal is to show that $H(x)$ has large oscillations.  Basically,
a zero of $\zeta(s)$ with real part $\b > 1/2$ induces oscillations
in $H(x)$ of size $x^{\b-1-\eps}$, which will overwhelm the error term
in \eqref{grbigosc}.

The Mellin transform of $H(x)$ does not exist because of the blow-up of
the integrand near $x=1$;  however the function $H(x)\log x$
is bounded near $x=1$.

\begin{lemma}\label{MellinH}
For $\Re s > 0$, we have
$$
\int_1^\infty x^{-s-1} H(x) \log x \, dx = - \frac{1}{s^2}
\log\pfrac{s \,
  \zeta(s+1)}{s+1} - \frac{1-\g}{s} + G(s),
$$
where $G(s)$ is a function that is analytic for $\Re s > -1$.
\end{lemma}

\begin{proof}  
By \eqref{logzeta},
$$
H(x)\log x = - (\log x) \sum_{n\le x} \frac{\Lambda(n)}{n\log n} + R(x) +
(\log x)(\log\log x + \g). 
$$
The Mellin transform of the sum is $s^{-1} \log \zeta(s+1)$, hence
$$
\int_1^\infty x^{-s-1} (\log x) \sum_{n\le x} \frac{\Lambda(n)}{n\log n}\, dx =
-\frac{d}{ds} \frac{\log \zeta(s+1)}{s} = \frac{\log \zeta(s+1)}{s^2} 
- \frac{1}{s} \frac{\zeta'}{\zeta}(s+1).
$$
Let
$$
f(x) \!:= \int_1^x \frac{1-t^{-1}}{t\log t}\, dt.
$$
We have (cf.~(6.7) of \cite{BD})
$$
f(x)\log x = (\log\log x + \g)\log x + O\pfrac{1}{x} \qquad (x > 1),
$$
and note that a piecewise continuous function which is $O(1/x)$ has a
Mellin transform which is analytic for $\Re s > -1$.  Also, 
$$
 s \int_1^\infty x^{-s-1} f(x)\, dx = \int_1^\infty x^{-s} f'(x)\, dx
= \int_1^\infty x^{-s} \frac{1-x^{-1}}{x\log x}\, dx = \log \pfrac{s+1}{s}.
$$
Thus, 
$$
\int_1^\infty x^{-s-1} f(x)\log x\, dx =
- \frac{d}{ds} \frac{1}{s} \log\pfrac{s+1}{s}  =
\frac{1}{s^2} \log\pfrac{s+1}{s} + \frac{1}{s^2} - \frac{1}{s} + \frac{1}{s+1}.
$$
Recalling \eqref{MTR}, the proof is complete.
\end{proof}

We see that the Mellin transform of 
$H(x)\log x$ has no real singularities in the
region $\Re s > -1$.  If $\zeta(s)$ has a zero with real part $\b>1/2$,
Landau's oscillation theorem implies that $H(x)\log x =
\Omega_\pm(x^{\b-1-\eps})$ for every $\eps>0$. 
Inequality \eqref{grbigosc} then implies that
$\g_r < \eg$ for infinitely many $r$, proving Theorem \ref{RHfalse}.

\begin{remark} We leave as an open problem to show that $\g_r > \eg$ for
infinitely many $r$ in case RH is false.  If the supremum $\sigma$ of
real parts of zeros of $\zeta(s)$ is strictly less than 1, then
Landau's oscillation theorem immediately gives
$$
H(x) = \Omega_\pm(x^{\sigma-1-\eps})
$$
for every $\eps>0$, while a simpler argument shows that
\[
R(x) = O(x^{\sigma-1+\eps}).
\]
By \eqref{gammaRHD}, we have
$$
\g_r = \eg \exp \left\{H(x) + O\pfrac{R^2(x)}{\log^2 x} + 
O\pfrac{1}{\sqrt{x}\log x} \right\}
$$
and the desired result follows immediately.  If $\zeta(s)$ has a sequence
of zeros with real parts approaching 1, Landau's theorem is too crude
to show that $H(x)$ has larger oscillations than does $R^2(x)/\log^2 x$.
In this case, techniques of Pintz (\cite{Pintz1}, \cite{Pintz2}) 
are perhaps useful.
\end{remark}


\emph{Acknowledement.}  The authors thank the referee for a careful reading
of the paper and for helpful suggestions regarding the exposition.

\end{document}